\input amstex
\magnification\magstephalf
\documentstyle{amsppt}


\hsize 5.72 truein
\vsize 7.9 truein
\hoffset .39 truein
\voffset .26 truein
\mathsurround 1.67pt
\parindent 20pt
\normalbaselineskip 13.8truept
\normalbaselines
\binoppenalty 10000
\relpenalty 10000
\csname nologo\endcsname 


\font\bc=cmb10
\font\tenbsy=cmbsy10

\catcode`\@=11

\def\qedsymbol{{\mathsurround\z@$\square$}}
\redefine\qed{\relaxnext@\ifmmode\let\next\@qed\else
  {\unskip\nobreak\hfil\penalty50\hskip2em\null\nobreak\hfil
    \qedsymbol\parfillskip\z@\finalhyphendemerits0\par}\fi\next}
\def\@qed#1$${\belowdisplayskip\z@\belowdisplayshortskip\z@
  \postdisplaypenalty\@M\relax#1
  $$\par{\lineskip\z@\baselineskip\z@\vbox to\z@{\vss\noindent\qed}}}
\outer\redefine\beginsection#1#2\par{\par\penalty-250\bigskip\vskip\parskip
  \leftline{\tenbsy x\bf#1. #2}\nobreak\smallskip\noindent}
\outer\redefine\genbeginsect#1\par{\par\penalty-250\bigskip\vskip\parskip
  \leftline{\bf#1}\nobreak\smallskip\noindent}

\def\next{\let\@sptoken= }\def\next@{ }\expandafter\next\next@
\def\@futureletnext#1{\let\nextii@#1\futurelet\next\@flti}
\def\@flti{\ifx\next\@sptoken\let\next@\@fltii\else\let\next@\nextii@\fi\next@}
\expandafter\def\expandafter\@fltii\next@{\futurelet\next\@flti}

\let\zeroindent\z@
\let\savedef@\endproclaim\let\endproclaim\relax 
\define\chkproclaim@{\add@missing\endroster\add@missing\enddefinition
  \add@missing\endproclaim
  \envir@stack\endproclaim
  \edef\endit@{\leftskip\the\leftskip\rightskip\the\rightskip}}
\let\endproclaim\savedef@
\def\thing@{.\enspace\egroup\ignorespaces}
\def\thingi@(#1){ \rm(#1)\thing@}
\def\thingii@\cite#1{ \rm\@pcite{#1}\thing@}
\def\thingiii@{\ifx\next(\let\next\thingi@
  \else\ifx\next\cite\let\next\thingii@\else\let\next\thing@\fi\fi\next}
\def\thing#1#2#3{\chkproclaim@
  \ifvmode \medbreak \else \par\nobreak\smallskip \fi
  \noindent\advance\leftskip#1
  \hskip-#1#3\bgroup\bc#2\unskip\@futureletnext\thingiii@}
\let\savedef@\endproclaim\let\endproclaim\relax 
\def\endit{\endproclaim\endit@\let\endit@\undefined}
\let\endproclaim\savedef@

\def\lemma#1{\thing\parindent{Lemma #1}\sl}

\def\thm#1{\thing\parindent{Theorem #1}\sl}

\def\remk#1{\thing\zeroindent{Remark #1}\rm}

\def\narrowthing#1{\chkproclaim@\medbreak\narrower\noindent
  \it\def\next{#1}\def\next@{}\ifx\next\next@\ignorespaces
  \else\bgroup\bc#1\unskip\let\next\narrowthing@\fi\next}
\def\narrowthing@{\@futureletnext\thingiii@}

\def\@cite#1,#2\end@{{\rm([\bf#1\rm],#2)}}
\def\cite#1{\in@,{#1}\ifin@\def\next{\@cite#1\end@}\else
  \relaxnext@{\rm[\bf#1\rm]}\fi\next}
\def\@pcite#1{\in@,{#1}\ifin@\def\next{\@cite#1\end@}\else
  \relaxnext@{\rm([\bf#1\rm])}\fi\next}

\def\pr@m@s{\ifx'\next\let\nxt\pr@@@s \else\ifx^\next\let\nxt\pr@@@t
  \else\let\nxt\egroup\fi\fi \nxt}

\input ifpdf.sty

\catcode`\@\active

\ifpdf
  \input supp-pdf
  \define\epsfbox#1{\hbox{\convertMPtoPDF{#1}{1}{1}}}
  \pdfpagewidth=8.5 truein
  \pdfpageheight=11 truein
\else
  \input epsf
\fi

\let\0\relax 
\define\codim{\operatorname{codim}}

\topmatter
\title On the Nochka-Chen-Ru-Wong Proof of Cartan's Conjecture\endtitle
\author Paul Vojta\endauthor
\affil University of California, Berkeley\endaffil
\address Department of Mathematics, University of California,
  970 Evans Hall\quad\#3840, Berkeley, CA \ 94720-3840;
  {\tt vojta\@math.berkeley.edu}\endaddress
\date 7 February 2008 \enddate
\thanks Supported by NSF grants DMS-9304899, DMS-0200892, and
  DMS-0500512.\endthanks
\subjclassyear{2000}
\subjclass Primary 11J25; Secondary 11J97, 32H30\endsubjclass
\keywords Nochka diagram; Schmidt's Subspace Theorem\endkeywords

\abstract In 1982--83, E. Nochka proved a conjecture of Cartan on defects
of holomorphic curves in $\Bbb P^n$ relative to a possibly degenerate set
of hyperplanes.  This was further explained by W. Chen in his 1987 thesis,
and subseqently simplified by M. Ru and P.-M. Wong in 1991.
The proof involved assigning weights to the hyperplanes.  This paper
provides further simplification of the proof of the construction of
the weights, by bringing back the use of the convex hull in working
with the ``Nochka diagram.''
\endabstract
\endtopmatter

\document
\beginsection{\01}{Introduction}

In 1982 and 1983, E. Nochka proved a conjecture of Cartan on defects of
holomorphic curves in $\Bbb P^n$ relative to a possibly degenerate set of
hyperplanes.  This was further explained by W. Chen in his thesis \cite{C},
and subsequently simplified by M. Ru and P.-M. Wong \cite{R-W}.
Ru and Wong also carried over the proof to the number field case,
giving an extension of Schmidt's Subspace Theorem.

In addition, Nochka recently published his original proof \cite{Noc}.
In that paper it is implicit but clear that he is using a convex hull
of a collection of points in the ``Nochka diagram,'' but this is not so
clear in \cite{R-W}.  On the other hand, the Ru-Wong proof avoids
Nochka's ``triangle inequalities'' and the definitions that they require.

This paper provides some further simplifications of the work of Nochka
and others, consisting of combining the simplifications of \cite{R-W}
with explicit use of a convex hull in defining the Nochka polygon,
as in \cite{V1}.  This proof also rewords the combinatorics to use
linear subspaces of $\Bbb P^n$ instead of sets of hyperplanes
(motivated by Shiffman's \cite{S} rephrasing of \cite{R-W} Thm.~2.2).

This paper only addresses the proof of the existence of the Nochka weights
(Theorem \01.1).  For details on the remainder of Nochka's proof, see
\cite{R-W}; simplified versions are also given in \cite{S} and \cite{V~2}.

Let $H_1,\dots,H_q$ be hyperplanes in $\Bbb P^k$, not necessarily distinct,
but in $n$\snug-subgeneral position; i.e., there exists an embedding of
$\Bbb P^k$ as a linear subspace of $\Bbb P^n$ and (distinct) hyperplanes
$H_1',\dots,H_q'$ in general position in $\Bbb P^n$ such that
$H_i=H_i'\cap\Bbb P^k$ for all $i$.

Nochka's theorem on the construction of Nochka weights is then:

\thm{\01.1} (Main Theorem)  If $q>2n-k+1$ then there exist weights
$\omega_1,\dots,\omega_q\in\Bbb R$ such that
\roster
\item"(i)."  $\omega_i\ge0$ for all $i$;
\item"(ii)."  $\omega_i\le\tau$ for all $i$, where
$$\tau = \frac{\sum_{i=1}^q\omega_i - k - 1}{q-2n+k-1};$$
and
\item"(iii)."  for any nonempty $L\subseteq\Bbb P^k$,
$$\sum_{\{i:H_i\supseteq L\}}\omega_i \le \codim L\;.\tag\01.1.1$$
\endroster
\endit

This then implies the following theorem \cite{R-W} Thm\. 3.5 (and its
counterpart in Nevanlinna theory).

\thm{\01.2}  Let $F$ be a number field, let $H_1,\dots,H_q$ be hyperplanes
in $\Bbb P^k_F$ in $n$\snug-subgeneral position (not necessarily distinct),
and let $\epsilon>0$.  Then there is a finite collection of proper linear
subspaces of $\Bbb P^k$ such that
$$\sum_{i=1}^q m(P,H_i) \le (2n-k+1+\epsilon)h(P)$$
for all $P\in\Bbb P^k(F)$ not lying in one of the linear subspaces.
\endit

\beginsection{\02}{The Nochka Diagram}

Before giving the proof of the Main Theorem, we first describe the
Nochka diagram of \cite{R-W}, and prove a lemma about the subspaces
occurring in it.

For linear subspaces $L\subseteq\Bbb P^k$, define
$$\alpha(L) = \#\{i:H_i\supseteq L\}$$
and recall that $\codim L$ denotes the codimension of $L$ in $\Bbb P^k$;
i.e., $\codim L=k-\dim L$.  Also, by convention, let $\codim\emptyset=k+1$.
For linear subspaces $L\subseteq\Bbb P^k$ let
$$P(L)=(\alpha(L),\codim L)\in\Bbb R^2\;.$$
By $n$\snug-subgeneral position, $P(L)$ lies above the line $\ell$ of slope $1$
passing through the point $(n,k)$ for all nonempty $L$; in other words,
$$\alpha(L) \le \codim L + n - k\tag\02.1$$
for all $L\ne\emptyset$.

Let $X$ be the point $(2n-k+1,k+1)$, and let $P_0,\dots,P_s,X$ be the lower
convex hull of the set
$$\{P(L):\emptyset\ne L\subseteq\Bbb P^k\}\cup\{X\}\;,$$
in order of increasing $x$\snug-coordinate.  For $j=0,\dots,s$ let $L_j$
be a linear subspace in $\Bbb P^k$ such that $P(L_j)=P_j$.  Let
$P_{s+1}=X$ and $L_{s+1}=\emptyset$, and note that $P(L_{s+1})\ne P_{s+1}$.
Also, $P_0=(0,0)$ and $L_0=\Bbb P^k$.

\centerline{\epsfbox{nochka.1}}

The following is motivated by \cite{R-W, Prop.~2.1}.
\lemma{\02.2}  For $j=0,\dots,s$, $L_j\supseteq L_{j+1}$.
\endit

\demo{Proof}  The case $j=s$ is trivial, so we may assume that $j<s$.

Note that $P_0,\dots,P_s$ must lie below the line $OX$ in the above
Nochka diagram; since they also must lie to the left of the line $\ell$,
they must therefore lie below and to the left of the point
$W=((2n-k+1)/2,(k+1)/2)$.  In particular, $\codim L_j\le(k+1)/2$.
(This is a determining factor in the placement of $X$.)
Also $\codim L_{j-1}<(k+1)/2$, so $\codim(L_{j-1}\cap L_j)<k+1$; in particular,
$L_{j-1}\cap L_j\ne\emptyset$.

We now use the facts that
$$\alpha(L_j+L_{j+1}) + \alpha(L_j\cap L_{j+1})
  \ge \alpha(L_j) + \alpha(L_{j+1})$$
and
$$\codim(L_j+L_{j+1}) + \codim(L_j\cap L_{j+1})
  = \codim L_j + \codim L_{j+1}\;.$$
It follows that at least one of the points $P(L_j+L_{j+1})$
or $P(L_j\cap L_{j+1})$ lies below or on the line $P_jP_{j+1}$.
By construction it therefore follows that $P(L_j+L_{j+1})=P(L_j)$
or $P(L_j\cap L_{j+1})=P(L_{j+1})$; hence $L_j+L_{j+1}=L_j$
or $L_j\cap L_{j+1}=L_{j+1}$.  Either of these conditions implies
that $L_j\supseteq L_{j+1}$.\qed
\enddemo

\beginsection{\03}{Proof of the Main Theorem}

We are now ready to prove the Main Theorem.

\demo{Proof of Main Theorem}
By Lemma \02.2 we may define $\omega_i$ to be the slope
of $P_{j-1}P_j$ for the smallest value of $j$ such that $H_i\supseteq L_j$.
Then condition (i) is trivially satisfied.

Let $\sigma$ be the slope of the line $P_sP_{s+1}$; then condition (ii)
is equivalent to the condition $\sigma\le\tau$.  But
$$\split \sum_{i=1}^q\omega_i
  &= \sum_{\{i:H_i\supseteq L_s\}}\omega_i
    + \sum_{\{i:H_i\nsupseteq L_s\}}\omega_i \\
  &= \codim L_s + \sigma(q-\alpha(L_s)) \\
  &= \codim L_s + \sigma(2n-k+1-\alpha(L_s)) + \sigma(q-2n+k-1) \\
  &= \codim L_s + (k+1-\codim L_s) + \sigma(q-2n+k-1) \\
  &= k+1+\sigma(q-2n+k-1).\endsplit$$
Thus, in fact, $\sigma=\tau$.

We now show condition (iii).  The proof is broken into two cases.

{\bc Case I:}  $L\cap L_s=\emptyset$.  In this case we have
$$\codim L+\codim L_s\ge k+1\;.\tag\03.1$$
By (\02.1) applied to $L_s$, (\03.1), and (\02.1) applied to $L$, we then have
$$\split \frac1\sigma &= \frac{2n-k+1-\alpha(L_s)}{k+1-\codim L_s} \\
  &\ge \frac{n+1-\codim L_s}{k+1-\codim L_s} \\
  &= 1 + \frac{n-k}{k+1-\codim L_s} \\
  &\ge 1 + \frac{n-k}{\codim L} \\
  &= \frac{\codim L + n - k}{\codim L} \\
  &\ge \frac{\alpha(L)}{\codim L}\;.\endsplit$$
Thus
$$\sum_{\{i:H_i\supseteq L\}}\omega_i \le \sigma\alpha(L) \le \codim L\;.$$

{\bc Case II:}  $L\cap L_s\ne\emptyset$.  We show, by induction on $j$,
that if $L\supseteq L_j$ then (\01.1.1) holds.  The case $j=s+1$ implies
condition (iii).  If $j=0$ then this claim is trivial.

Suppose now that $L\supseteq L_j$ and that (\01.1.1) holds for $L+L_{j-1}$:
$$\sum_{\{i:H_i\supseteq L+L_{j-1}\}}\omega_i \le \codim(L+L_{j-1})\;.
  \tag\03.2$$
Let $\sigma_{j-1}$ denote the slope of the line $P_{j-1}P_j$.
Since $L\cap L_{j-1}\supseteq L\cap L_s\ne\emptyset$ by assumption,
the point $P(L\cap L_{j-1})$ lies on or above the line $P_{j-1}P_j$.  Thus
$$\split \codim(L\cap L_{j-1}) - \codim L_{j-1}
  &\ge \sigma_{j-1}(\alpha(L\cap L_{j-1})-\alpha(L_{j-1})) \\
  &=\sum_{\{i:\text{$H_i\supseteq L\cap L_{j-1}$ and $H_i\nsupseteq L_{j-1}$}\}}
    \omega_i \\
  &= \sum_{\{i:H_i\supseteq L\cap L_{j-1}\}}\omega_i
    - \sum_{\{i:H_i\supseteq L_{j-1}\}}\omega_i\;,
\endsplit\tag\03.3$$
where the second step uses the assumption that $L\supseteq L_j$.

We also have
$$\sum_{\{i:H_i\supseteq L\}}\omega_i
    + \sum_{\{i:H_i\supseteq L_{j-1}\}}\omega_i
  \le \sum_{\{i:H_i\supseteq L+L_{j-1}\}}\omega_i
    + \sum_{\{i:H_i\supseteq L\cap L_{j-1}\}}\omega_i$$
and
$$\codim L + \codim L_{j-1}
  = \codim(L + L_{j-1}) + \codim(L\cap L_{j-1})\;.$$
Therefore (\03.2) and (\03.3) give
$$\split \codim L &= \codim(L+L_{j-1}) + \codim(L\cap L_{j-1})
    - \codim L_{j-1} \\
  &\ge \sum_{\{i:H_i\supseteq L+L_{j-1}\}}\omega_i
    + \sum_{\{i:H_i\supseteq L\cap L_{j-1}\}}\omega_i
    - \sum_{\{i:H_i\supseteq L_{j-1}\}}\omega_i \\
  &\ge \sum_{\{i:H_i\supseteq L\}}\omega_i\;.\endsplit$$
This gives (\01.1.1) and therefore the theorem is proved.\qed
\enddemo

\remk{\03.4}  Theorem \01.1 usually has an additional condition
$\tau\le(k+1)/(n+1)$.  This is easy to see from the Nochka diagram, since
the line from the point $(n-k,0)$ (where $\ell$ meets the $x$\snug-axis)
to $X$ has slope $(k+1)/(n+1)$.  By a careful examination of Nochka's proof,
however, N. Toda \cite{Nog, p.~340} has noted that this can be improved to
$\tau\le k/n$.  To see this using the present proof, let $(x,y)$ be the
coordinates of $P_s$; then Toda's condition is equivalent to
$$\frac{k+1-y}{2n-k+1-x} \le \frac kn\;.\tag\03.4.1$$
By (\02.1), we have
$$kx-ky \le k(n-k)\;,$$
and the fact that $y\ge1$ gives
$$-(n-k)y \le -(n-k)\;.$$
Adding these two inequalities then gives
$$kx-ny \le (k-1)(n-k)\;,$$
which is equivalent to (\03.4.1).
\endit

\enddocument